\documentclass[12pt,leqno]{article}
\usepackage{amssymb,eufrak}
\begin{document}
\centerline{\bf On the $p$-adic Leopoldt transform of a power series}

\centerline{ Bruno Angl\`  es\footnote{Universit\'e de Caen, 
LMNO  CNRS UMR 6139, 
BP 5186, 14032 Caen Cedex, France. E-mail: angles@math.unicaen.fr} }
${}$\par
%\noindent Universit\'e de Caen, 
%Laboratoire de Math\'ematiques Nicolas Oresme, CNRS UMR 6139, 
%Campus II, Boulevard Mar\'echal Juin,
%BP 5186, 14032 Caen Cedex, France.\par
%\noindent E-mail: angles@math.unicaen.fr
%\date{April 2007}

%\begin{abstract}
 %\end{abstract}\par
 %\noindent MSC: 11R23; 11R18\par
 %\noindent{\sl  Keywords: $p$-adic L-functions; Iwasawa theory; Cyclotomic fields}\par
 %${}$\par

 %-------------------------------------INTRODUCTION--------------------------------------------------
 Let $p$ be an odd prime number. Let $X$ be the projective limit  for the norm maps of the $p$-Sylow subgroups of the ideal class groups of $\mathbb Q(\zeta_{p^{n+1}}),$ $n\geq 0.$ Let $\Delta ={\rm Gal}(\mathbb Q(\zeta_p)/\mathbb Q)$ and let $\theta $ be an even and non-trivial character of $\Delta. $ Then $X$ is a $\mathbb Z_p[[T]]$-module and the characteristic ideal of the isotypic component $X(\omega \theta^{-1})$ is generated by a power series $f(T,\theta)\in \mathbb Z_p[[T]]$ such that (see for example \cite{CS}):
 $$\forall n\geq 1,\, n\equiv 0\pmod{p-1},\, f((1+p)^{1-n}-1,\theta)=L(1-n,\theta ),$$
 where $L(s,\theta)$ is the usual Dirichlet $L$-series.  Therefore, it is natural and interesting to study the properties of the power series $f(T,\theta).$\par
 ${}$\par
 We denote by $\overline{f(T,\theta)}\in \mathbb F_p[[T]]$ the reduction of $f(T,\theta)$ modulo $p.$ Then B. Ferrero and L. Washington have proved (\cite{FW}):
 $$\overline{f(T,\theta)}\not =0.$$
 Note that, in fact, we have (\cite{ANG}):
 $$\overline{f(T,\theta)}\not \in \mathbb F_p[[T^p]].$$
W. Sinnott has  proved the following (\cite{SI2}):
 $$\overline{f(T,\theta)}\not \in \mathbb F_p(T).$$
 But, note that $\forall a\in \mathbb Z_p^*,$  $\mathbb F_p[[T]]=\mathbb F_p[[(1+T)^a -1]].$ Therefore it is natural to introduce the notion of a pseudo-polynomial which is an element  $F(T)$ in $\mathbb F_p[[T]]$ such that there exist an integer $r\geq 1,$ $c_1,\cdots c_r\in \mathbb F_p,$ $a_1,\cdots ,a_r\in \mathbb Z_p,$ such that $F(T)=\sum_{i=1}^r c_i (1+T)^{a_i}.$ An element of $\mathbb F_p[[T]]$ will be called a pseudo-rational function if it is the quotient of two pseudo-polynomials. In this paper, we prove that $\overline{f(T,\theta)}$ is not a pseudo-rational function (part 1) of Theorem \ref{Theorem2}). This latter result suggests the following question: is $\overline{f(T,\theta)}$ algebraic over $\mathbb F_p(T)?$ We suspect that  this  is not the case but we have no evidence for it. Note that, by the result of Ferrero and Washington, we can write:
 $$\overline{f(T,\theta)}=T^{\lambda (\theta )}U(T),$$
 where $\lambda (\theta) \in \mathbb N $ and $U(T)\in \mathbb F_p[[T]]^*.$ S. Rosenberg has proved that (\cite{ROS}):
 $$\lambda (\theta )\leq (4p(p-1))^{\phi (p-1)},$$
 where $\phi$ is Euler's totient function. In this paper, we improve Rosenberg's bound  (part 2) of Theorem \ref{Theorem2}):
 $$\lambda(\theta )<(\frac{p-1}{2})^{\phi (p-1)}.$$
 This implies that the lambda invariant of the field $\mathbb Q(\zeta_p)$ is less than $2(\frac{p-1}{2})^{\phi (p-1)+1}$ (see Corollary \ref{Theorem3} for the precise statement for an abelian number field). Note that this bound is certainly far from the truth, because according to a heuristic argument due to Ferrero and Washington (see \cite{LA})  and to Grennberg's conjecture:
 $$\lambda (\mathbb Q(\zeta_p))=\sum_{\theta \in \widehat{\Delta},\, \theta \not =1\, {\rm  and \, even }}\lambda(\theta )\leq \frac{{\rm Log}(p)}{{\rm Log}({\rm Log}(p)) }.$$
 ${}$\par
 The author is indebted to Warren Sinnott for communicating some of his unpublished works (note that Lemma \ref{Lemma7} is due to Warren Sinnott). The author also thanks Filippo Nuccio for pointing out the work of J. Kraft and L. Washington (\cite{KW}).\par

 %--------------------------------------NOTATIONS--------------------------------------------------------

 \section{Notations}\par
 ${}$\par
 Let $p$ be an odd prime number and let $K$ be a finite extension of $\mathbb Q_p.$ Let $O_K$ be the valuation ring of $K$ and let $\pi $ be a prime of $K.$ We set $\mathbb F_q=O_K/\pi O_K,$  it is a finite field having $q$ elements and its characteristic is $p.$ Let $T$ be an indeterminate over $K,$ we set $\Lambda =O_K[[T]].$ Observe that $\Lambda/\pi \Lambda \simeq \mathbb F_q[[T]].$ Let $F(T)\in \Lambda\setminus \{ 0\},$ then we can write in an unique way (\cite{WAS}, Theorem 7.3):
 $$F(T)=\pi^{\mu(F)} P(T) U(T),$$
 where $U(T)$ in an unit of $\Lambda ,$ $\mu(F)\in \mathbb N,$ $P(T)\in O_K[T]$ is a monic polynomial such that $P(T)\equiv T^{\lambda (F)}\pmod{\pi }$ for some integer $\lambda (F)\in \mathbb N.$ If $F(T)=0,$ we set $\mu(F)=\lambda(F)=\infty .$ An element $F(T)\in \Lambda$ is called a pseudo-polynomial (see also \cite{ROS}, Definition 2) if there exist some integer $r\geq 1,$ $c_1,\cdots , c_r \in O_K,$ $a_1,\cdots ,a_r\in \mathbb Z_p,$ such that:
 $$F(T)=\sum_{i=1}^r c_i (1+T)^{a_i}.$$
 We denote the ring of pseudo-polynomials in $\Lambda $ by $A.$ Let $\delta \in \mathbb Z/(p-1)\mathbb Z$ and $F(T)\in \Lambda ,$ we set:
 $$\gamma_{\delta } (F(T))= \frac{1}{p-1} \sum_{\eta \in \mu_{p-1}} \eta^{\delta } F((1+T)^{\eta }-1).$$
 Then $\gamma_{\delta}:\Lambda \rightarrow \Lambda$ is a $O_K$-linear map and:\par
 \noindent - for $\delta ,\delta'\in \mathbb Z/(p-1)\mathbb Z,$ $\gamma_{\delta}\gamma_{\delta'}=0$ if $\delta \not = \delta'$ and $\gamma_{\delta}^2=\gamma_{\delta},$\par
 \noindent - $\sum_{\delta \in \mathbb Z/(p-1)\mathbb Z} \gamma_{\delta }={\rm Id}_{\Lambda}.$\par
 \noindent  For $F(T)\in \Lambda,$ we set:
 $$D(F(T))=(1+T) \frac{d}{dT} F(T),$$
 $$U(F(T))=F(T)-\frac{1}{p}\sum_{\zeta\in \mu_p} F(\zeta (1+T)-1)\, \in \Lambda .$$
  Then $D,U: \Lambda \rightarrow \Lambda $ are $O_K$-linear maps. Observe that:\par 
    \noindent - $U^2=U,$\par
  \noindent - $DU=UD,$\par
  \noindent- $\forall \delta \in \mathbb Z/(p-1)\mathbb Z,$ $ \gamma_{\delta }U=U\gamma_{\delta},$\par
  \noindent - $\forall \delta \in \mathbb Z/(p-1)\mathbb Z,$ $D\gamma_{\delta }=\gamma_{\delta +1}D.$\par
  \noindent If $F(T)\in \Lambda,$ we denote its reduction modulo $\pi$ by $\overline{F(T)}\in \mathbb F_q[[T]].$ If $f:\Lambda \rightarrow \Lambda$ is a $O_K$-linear map, we denote its reduction modulo $\pi$ by $\overline{f}: \mathbb F_q[[T]]\rightarrow \mathbb F_q[[T]].$ For all $n\geq 0,$ we set $\omega_n(T)=(1+T)^{p^n}-1.$\par
  Let $B$ be a commutative and unitary ring. We denote the set of invertible elements of $B$ by $B^*.$\par
   We fix $\kappa $ a topological generator of $1+p\mathbb Z_p.$ Let $x\in \mathbb Z_p$ and let $n\geq 1,$ we denote the unique integer $k\in \{ 0,\cdots , p^n-1\}$ such that $x\equiv k\pmod{p^n}$ by $[x]_n.$ Let $\omega :\mathbb Z_p^*\rightarrow \mu_{p-1}$ be the Teichm\"uller character, i.e. $\forall a\in \mathbb Z_p^*,$ $\omega (a)\equiv a \pmod{p}.$ Let $x,y\in \mathbb Z_p,$ we write:\par
   \noindent - $x\sim y$ if there exists $\eta \in \mu_{p-1}$ such that $y=\eta x,$\par
   \noindent - $x\equiv y \pmod {\mathbb Q^*}$ if there exists $z\in \mathbb Q^*$ such that $y=zx.$\par
   \noindent The function ${\rm Log}_p$ will denote the usual $p$-adic logarithm.  $v_p$ will denote the usual $p$-adic valuation on $\mathbb C_p$ such that $v_p(p)=1. $\par
   Let $\rho$ be a Dirichlet character of conductor $f_{\rho}.$ Recall that the Bernoulli numbers $B_{n,\rho}$ are defined by the following identity:
   $$\sum_{a=1}^{f_{\rho}}\frac{\rho (a) e^{aZ}}{e^{fZ}-1}=\sum_{n\geq 0} \frac{B_{n,\rho }}{n!} Z^{n-1},$$
   where $e^Z=\sum_{n\geq 0} Z^n/n!.$  If $\rho =1,$ for $n\geq 2,$ $B_{n,1}$ is the $n$th Bernoulli number.\par
   Let $x\in \mathbb R.$ We denote the biggest integer less than or equal to $x$ by $[x].$ The function ${\rm Log}$ will denote the usual logarithm.\par

 %--------------------------------------PRELIMINARIES----------------------------------------------------

 \section{Preliminaries}
 ${}$\par

 Let $\delta \in \mathbb Z/(p-1)\mathbb Z.$ In this section, we will recall the construction of the $p$-adic Leopoldt transform $\Gamma_{\delta }$ (see \cite{LA}, Theorem 6.2) which is a $O_K$-linear map from $\Lambda $ to $\Lambda .$\par
 First, observe that $(\pi^n, \omega_n(T))= \pi^n \Lambda +\omega_n(T)\Lambda ,$ $n\geq 1,$ is a basis of neighbourhood of zero in $\Lambda :$\par
 \newtheorem{Lemma1}{Lemma}[section]
 \begin{Lemma1} \label{Lemma1}
 ${}$\par
 \noindent 1) $\forall n\geq 1,$ $(\pi , T)^{2n}\subset (\pi^n, T^n)\subset (\pi ,T)^n.$\par
 \noindent 2) $\forall n\geq 1,$ $\omega_n (T)\in (p^{[n/2]}, T^{p^{[n/2]+1}}).$\par
 \noindent 3) Let $N\geq 1,$ set $n=[{\rm Log}(N)/{\rm Log}(p)].$ We have:
 $$T^N\in (p^{[n/2]}, \omega_{[n/2]+1}(T)).$$
 \end{Lemma1}
 \noindent{\sl Proof}  Note that assertion 1) is obvious. Assertion 2) comes from the fact:
 $$\forall k\in \{ 1,\cdots, p^n\},\, v_p(   \frac{p^n!}{k! (p^n-k)!})=n-v_p(k).$$
 To prove assertion 3), it is enough to prove the following:\par
 \noindent $\forall n\geq 0,$ there exist $\delta_0^{(n)}(T),\cdots , \delta_n^{(n)}(T) \in \mathbb Z [T]$ such that:
 $$T^{p^n}=\sum_{i+j=n} \omega_i(T) p^j \delta_j^{(n)}(T).$$
  Let's prove this latter fact by recurrence on $n.$ Note that the result is clear if $n=0.$ Let's assume that it is true for $n$ and let's prove the assertion for $n+1.$ Let $r(T)\in \mathbb Z[T]$ such that:
  $$\frac{\omega_{n+1}(T)}{\omega_n(T)} +pr(T)=T^{p^n(p-1)}.$$
  Then:
  $$T^{p^{n+1}}=T^{p^n} \frac{\omega_{n+1}(T)}{\omega_n(T)} +pr(T) T^{p^n}.$$
 Note that there exists $q(T) \in \mathbb Z[T]$ such that:
 $$\frac{\omega_{n+1}(T)}{\omega_n(T)}=\omega_n(T)^{p-1}+pq(T).$$
 Thus:
 $$T^{p^{n+1}}=\omega_{n+1}(T)\delta_0^{(n)}(T) +\, \sum_{i+j=n,\, j\geq 1}( \omega_n(T)^{p-1}+pq(T)) \omega_i(T) p^j \delta_j^{(n)}(T)\, + \sum_{i+j=n} \omega_i(T) p^{j+1} \delta_j^{(n)}(T) r(T).$$
 Thus, there exist $\delta_0^{(n+1)}(T), \cdots , \delta_{n+1}^{(n+1)}(T) \in \mathbb Z[T]$ such that:
 $$T^{p^{n+1}}=\sum_{i+j=n+1} \omega_i(T)p^j \delta_j^{(n+1)} (T).\, \diamondsuit $$
 The following Lemma will be useful in the sequel (for a similar result see \cite{ROS}, Lemma 5):\par
  \newtheorem{Lemma2}[Lemma1]{Lemma}
 \begin{Lemma2} \label{Lemma2}
 Let $F(T)\in A.$ Write $F(T)=\sum_{i=1}^r \beta_i (1+T)^{\alpha_i},$ $\beta_1,\cdots, \beta_r\in O_K,$ $\alpha_1,\cdots ,\alpha_r \in \mathbb Z_p,$ and $\alpha_i\not = \alpha_j$ for $i\not = j.$ Let $N={\rm Max}\{ v_p(\alpha_i-\alpha_j),i\not =j\}.$  Let $n\geq 1$ be an integer. Then:
 $$F(T)\equiv 0\pmod{(\pi^n,\omega_{N+1}(T))} \Leftrightarrow \forall i=1,\cdots r, \, \beta_i\equiv 0\pmod{\pi^n}.$$
 \end{Lemma2}
 \noindent{\sl Proof} We have:
 $$F(T)\equiv \sum_{i=1}^r \beta_i (1+T)^{[\alpha_i]_{N+1}}\pmod{\omega_{N+1}(T)}.$$
 Therefore $F(T)\equiv 0\pmod{(\pi^n,\omega_{N+1}(T))}$ if and only if we have:
 $$\sum_{i=1}^r\beta_i (1+T)^{[\alpha_i]_{N+1}}\equiv 0\pmod{\pi^n}.$$
 But for $i\not =j,$ $[\alpha_i]_{N+1}\not = [\alpha_j]_{N+1}.$ Therefore $\sum_{i=1}^r\beta_i (1+T)^{[\alpha_i]_{N+1}}\equiv 0\pmod{\pi^n}$ if and only if:
 $$\forall i=1,\cdots r, \, \beta_i\equiv 0\pmod{\pi^n}.\, \diamondsuit $$
 Observe that $U,D, \gamma_{\delta}$ are continuous $O_K$-linear maps by Lemma \ref{Lemma1} and the following Lemma:\par
 \newtheorem{Lemma3}[Lemma1]{Lemma}
 \begin{Lemma3} \label{Lemma3}
Let $F(T)\in \Lambda$ and let $n\geq 0.$\par
\noindent 1) $F(T)\equiv 0\pmod{\omega_n(T)}\Rightarrow \gamma_{\delta}(F(T))\equiv 0\pmod{\omega_n(T)}.$\par
\noindent 2) $F(T)\equiv 0\pmod{\omega_n(T)}\Rightarrow D(F(T))\equiv 0\pmod{(p^n,\omega_n(T))}.$\par
\noindent 3) If $n\geq 1,$ $F(T)\equiv 0\pmod{\omega_n(T)}\Rightarrow U(F(T))\equiv 0\pmod{\omega_n(T)}.$\par
\end{Lemma3}
\noindent{\sl Proof} The assertions 1) and 2) are obvious. It remains to prove 3). Observe that, by \cite{WAS}, Proposition 7.2, we have:
$$\forall G(T)\in \Lambda,\, G(T)\equiv 0\pmod{\omega_n(T)} \Leftrightarrow \forall \zeta \in \mu_{p^n},\, G(\zeta-1)=0.$$
 Now, let $F(T)\in \Lambda ,$ $F(T)\equiv 0\pmod{\omega_n(T)}.$ Since the map: $\mu_{p^n}\rightarrow \mu_{p^n},$ $x\mapsto \zeta x,$ is a bijection for all $\zeta\in \mu_{p^n},$ we get:
 $$\forall \zeta\in \mu_{p^n},\, U(F)(\zeta-1)=0.$$
 Therefore: 
 $$U(F(T))\equiv 0\pmod{\omega_n(T)}.\, \diamondsuit $$
 Let $s\in \mathbb Z_p.$ For $n\geq 0,$ set:
 $$k_n(s,\delta)=[s]_{n+1}+\delta_n p^{n+1}\, \in \mathbb N\setminus\{ 0\},$$
 where $\delta_n\in \{ 1, \cdots ,p-1\}$ is such that $[s]_{n+1}+\delta_n \equiv \delta \pmod{p-1}.$ Observe that:\par
 \noindent - $\forall n\geq 0,$ $k_n(s,\delta)\equiv \delta \pmod{p-1}$ and $k_n(s,\delta)\equiv s\pmod{p^{n+1}},$\par
 \noindent - $\forall n\geq 0,$ $k_{n+1}(s,\delta)>k_n(s,\delta),$\par
 \noindent - $s={\rm lim}_nk_n(s,\delta).$\par
 \noindent  In particular:
 $$\forall a\in \mathbb Z_p,\, \forall n\geq 0,\, a^{k_{n+1}(s,\delta)}\equiv a^{k_n(s,\delta)}\pmod{p^{n+1}}.$$
 Now, let $F(T)\in A.$ Write $F(T)=\sum_{i=1}^r \beta_i (1+T)^{\alpha_i},$ $\beta_1,\cdots ,\beta_r\in O_K,$ $\alpha_1,\cdots ,\alpha_r \in \mathbb Z_p.$ We set:
 $$\Gamma_{\delta }(F(T))=\sum_{\alpha_i\in \mathbb Z_p^*} \beta_i \omega^{\delta}(\alpha_i) (1+T)^{\frac{{\rm Log}_p(\alpha_i)}{{\rm Log}_p(\kappa)}}.$$
 Thus, we have a surjective $O_K$-linear map: $\Gamma_{\delta }: A\rightarrow A.$ Note that:\par
 \newtheorem{Lemma4}[Lemma1]{Lemma}
 \begin{Lemma4} \label{Lemma4}
 Let $F(T)\in A.$\par
 \noindent 1) Let $s\in \mathbb Z_p.$ Then:
 $$\forall n\geq 0, \, \Gamma_{\delta} (F)(\kappa^s-1)\equiv D^{k_n(s,\delta)}(F)(0)\pmod{p^{n+1}}.$$
 2) Let $n\geq 1.$ Assume that $F(T)\equiv 0\pmod{\omega_n(T)}.$ Then $\Gamma_{\delta }(F(T))\equiv 0\pmod{\omega_{n-1}(T)}.$\par
 \end{Lemma4}
 \noindent {\sl Proof}  For $a\in \mathbb Z_p^*,$ write $a=\omega (a)\, <a>,$ where $<a>\in 1+p\mathbb Z_p.$ Let's write:
 $$F(T)=\sum_{i=1}^r\beta_i(1+T)^{\alpha_i},$$
 $\beta_1,\cdots ,\beta_r \in O_K,$ $\alpha_1,\cdots ,\alpha_r \in \mathbb Z_p.$ We have:
 $$D^{k_n(s,\delta)}(F(T))=\sum_{i=1}^r\beta_i \alpha_i^{k_n(s,\delta)}(1+T)^{\alpha_i}.$$
 Thus:
 $$D^{k_n(s,\delta)}(F(T))\equiv \sum_{\alpha_i \in \mathbb Z_p^*}\beta_i \omega^{\delta }(\alpha_i)<\alpha_i>^s(1+T)^{\alpha_i}\pmod{p^{n+1}}.$$
 But recall that:
 $$\Gamma_{\delta}(F)(\kappa^s-1)=\sum_{\alpha_i\in \mathbb Z_p^*}\beta_i \omega^{\delta}(\alpha_i)<\alpha_i>^s.$$
 Assertion 1) follows easily.  Now, let's suppose that $F(T)\equiv 0\pmod{\omega_n(T)}$ for some $n\geq 1.$ Then:
 $$\forall a\in\{ 0,\cdots ,p^n-1\},\, \sum_{\alpha_i\equiv a\pmod{p^n}}\beta_i =0.$$
 This implies that:
 $$\forall a\in\{ 0,\cdots, p^{n-1}-1\},\, \sum_{\alpha_i\in \mathbb Z_p^*,\, {\rm Log}_p(\alpha_i)/{\rm Log}_p(\kappa)\equiv a \pmod{p^{n-1}}}\omega^{\delta }(\alpha_i) \beta_i =0.$$
 But recall that:
 $$\Gamma_{\delta }(F(T))=\sum_{\alpha_i\in \mathbb Z_p^*} \beta_i \omega^{\delta}(\alpha_i) (1+T)^{\frac{{\rm Log}_p(\alpha_i)}{{\rm Log}_p(\kappa)}}.$$
 Thus  $\Gamma_{\delta }(F(T))\equiv 0\pmod{\omega_{n-1}(T)}.\, \diamondsuit $\par
 \newtheorem{Proposition1}[Lemma1]{Proposition}
 \begin{Proposition1} \label{Proposition1}
 Let $F(T)\in \Lambda.$ There exists an unique power series $\Gamma_{\delta}(F(T))\in \Lambda$ such that:
 $$\forall s\in \mathbb Z_p\, \forall n\geq 0,\, \Gamma_{\delta} (F)(\kappa^s-1)\equiv D^{k_n(s,\delta)}(F)(0)\pmod{p^{n+1}}.$$
 \end{Proposition1}
 \noindent{\sl Proof} Let $(F_N(T))_{N\geq 0}$ be a sequence of elements in $A$ such that:
 $$\forall N\geq 0,\, F(T)\equiv F_N(T)\pmod{\omega_N(T)}.$$
 Fix $N\geq 1.$ Then:
 $$\forall m\geq N, F_m(T)\equiv F_N(T)\pmod{\omega_N(T)}.$$
 Therefore, by Lemma \ref{Lemma4}, we have:
 $$\forall m\geq N, \Gamma_{\delta}(F_m(T))\equiv \Gamma_{\delta}(F_N(T))\pmod{\omega_{N-1}(T)}.$$
 This implies that the sequence $(\Gamma_{\delta}(F_N(T)))_{N\geq 1}$ converges in $\Lambda $ to some power series $G(T)\in \Lambda.$ Observe that, since $\Lambda $ is compact, we have:
 $$\forall N\geq 1, G(T)\equiv \Gamma_{\delta}(F_N(T))\pmod{\omega_{N-1}(T)}.$$
 In particular:
 $$\forall N\geq 1, G(\kappa^s-1)\equiv \Gamma_{\delta }(F_N)(\kappa^s-1)\pmod{p^N}.$$
 Thus, applying Lemma \ref{Lemma4}, we get:
 $$\forall N\geq 1, G(\kappa^s-1)\equiv D^{k_{N-1}(s,\delta )}(F_N)(0)\pmod{p^N}.$$
 But:
 $$\forall N\geq 1, D^{k_{N-1}(s,\delta )}(F(T))\equiv D^{k_{N-1}(s,\delta )}(F_N(T))\pmod{(p^N,\omega_N(T))}.$$
 Therfore:
 $$\forall N\geq 1, G(\kappa^s -1)\equiv D^{k_{N-1}(s,\delta )}(F)(0)\pmod{p^N}.$$
 Now, set $\Gamma_{\delta }(F(T))=G(T).$ The Proposition follows easily. $\diamondsuit$\par

 %--------------------------------------SOME PROPERTIES-----------------------------------------------

 \section{Some properties of the $p$-adic Leopoldt transform}
 ${}$\par
 We need the following fundamental result:\par
 \newtheorem{Proposition2}{Proposition}[section]
 \begin{Proposition2} \label{Proposition2}
 Let $\delta \in \mathbb Z/(p-1)\mathbb Z$ and let $F(T)\in \Lambda.$ Let $m,n\in \mathbb N\setminus \{ 0\}.$ then:
 $$\Gamma_{\delta}(F(T))\equiv 0\pmod{(\pi^n,\omega_{m-1}(T))} \Leftrightarrow \gamma_{-\delta}U(F(T)) \equiv 0\pmod{(\pi^n,\omega_{m}(T))}.$$
 \end{Proposition2}
 \noindent{\sl Proof}  A similar result has been obtained by S. Rosenberg (\cite{ROS}, Lemma 8). We begin by proving that $\Gamma_{\delta}$ is a continuous $O_K$-linear map. By Lemma \ref{Lemma1},  this comes from the following fact:\par
 \noindent Let $F(T)\in \Lambda.$ Let $n\geq 1$ and assume that $F(T)\equiv 0 \pmod{\omega_n(T)},$ then $\Gamma_{\delta }(F(T))\equiv 0\pmod{\omega_{n-1}(T)}.$\par
 \noindent Indeed, let $(F_N(T))_{N\geq 0}$ be a sequence of elements in $A$ such that:
 $$\forall N\geq 0, \, F(T)\equiv F_N(T) \pmod{\omega_{N}(T)}.$$
 By the proof of Proposition \ref{Proposition1}:
 $$\forall N\geq 1,\, \Gamma_{\delta}(F(T))\equiv \Gamma_{\delta}(F_N(T))\pmod{\omega_{N-1}(T)}.$$
 Now, by Lemma \ref{Lemma4}:
 $$\Gamma_{\delta}(F_n(T))\equiv 0\pmod{\omega_{n-1}(T)}.$$
 The assertion follows.\par
 \noindent  Now, since $\Gamma_{\delta}, \gamma_{-\delta},U$ are continuous $O_K$-linear maps, it suffices to prove the Proposition in the case where $F(T)\in A.$ Write $F(T)=\sum_{i=1}^r \beta_i (1+T)^{\alpha_i},$ $\beta_1,\cdots ,\beta_r \in O_K,$ $\alpha_1,\cdots ,\alpha_r \in \mathbb Z_p.$ Let $I\subset \{ \alpha_1,\cdots ,\alpha_r\}$ be a set of representatives of the classes of $\alpha_1, \cdots ,\alpha_r$ for the relation $\sim .$ For $x\in I,$ $x\not \equiv 0\pmod{p},$ set:
 $$\beta_x=\sum_{\alpha_i\sim x}\beta_i \frac{\alpha_i}{x}.$$
 We get:
 $$(p-1)\gamma_{-\delta}U(F(T))=\sum_{\eta\in \mu_{p-1}}\sum_{x\in I,\, x\in \mathbb Z_p^*}\eta^{-\delta}\beta_x (1+T)^{\eta x}.$$
 Now observe that:
 $$\Gamma_{\delta}(F(T)) =\Gamma_{\delta} \gamma_{-\delta }U (F(T))=\sum_{x\in I, \, x\in \mathbb Z_p^*}\beta_x \omega^{\delta}(x) (1+T)^{{\rm Log}_p(x)/{\rm Log}_p(\kappa )}.$$
 Therefore  $\Gamma_{\delta}(F(T))\equiv 0\pmod{(\pi^n,\omega_{m-1}(T))} $ if and only if:
 $$\forall a \in \{ 0,\cdots p^{m-1}-1\},\, \sum_{x\in I,\, x\in \mathbb Z_p^*,\, {\rm Log}_p(x)/{\rm Log}_p(\kappa )\equiv a\pmod{p^{m-1}}}\beta_x \omega^{\delta} (x) \equiv 0\pmod{\pi^n}.$$
 Now, observe that for $a\in \{ 0,\cdots , p^{m}-1\},$ there exists at most one $\eta \in \mu_{p-1}$ such that $[\eta x]_m=a,$ and if such a $\eta $ exists it is equal to $\omega (a) \omega^{-1}(x).$ Therefore  $\Gamma_{\delta}(F(T))\equiv 0\pmod{(\pi^n,\omega_{m-1}(T))} $ if and only if:
 $$\forall a \in \{ 0,\cdots ,p^m-1\},\, \sum _{x\in I,\, x\in \mathbb Z_p^*,\, \exists \eta_x \in \mu_{p-1},[\eta_x x]_m=a} \beta_x \eta_x^{-\delta}\equiv 0\pmod{\pi ^n}.$$
 This latter property is equivalent to $ \gamma_{-\delta}U(F(T)) \equiv 0\pmod{(\pi^n,\omega_{m}(T))}.\, \diamondsuit $\par
 Now, we can list the basic properties of $\Gamma_{\delta}:$\par
 \newtheorem{Proposition3}[Proposition2]{Proposition}
 \begin {Proposition3} \label{Proposition3}
 Let $\delta \in \mathbb Z/(p-1)\mathbb Z.$\par
 \noindent 1) $\Gamma_{\delta}:\Lambda \rightarrow \Lambda$ is a surjective and continuous $O_K$-linear map.\par
 \noindent 2) $\forall F(T)\in \Lambda,$ $\Gamma_{\delta}(F(T))=\Gamma_{\delta }\gamma_{-\delta}U(F(T)).$\par
 \noindent 3) $\forall a \in \mathbb Z_p^*,$ $\Gamma_{\delta} (F((1+T)^a-1))= \omega^{\delta }(a) (1+T)^{{\rm Log}_p(a)/{\rm Log}_p (\kappa)}\Gamma_{\delta }(F(T)).$\par
 \noindent 4) Let $\kappa'$ be another topological generator of $1+p\mathbb Z_p$ and let $\Gamma_{\delta}'$ be the $p$-adic Leopoldt transform associated to $\kappa'$ and $\delta .$ Then:
 $$\forall F(T)\in \Lambda ,\,  \Gamma_{\delta}'(F(T)) =\Gamma_{\delta}(F)((1+T)^{{\rm Log}_p(\kappa)/{\rm Log}_p(\kappa')}-1).$$
 5) Let $F(T)\in \Lambda.$ Then $\mu(\Gamma_{\delta }(F(T)))=\mu (\gamma_{-\delta}U(F(T)))$ and:
 $$\forall N\geq 1,\, \lambda (\Gamma_{\delta}(F(T)))\geq p^{N-1} \Leftrightarrow \lambda (\gamma_{-\delta}U(F(T)))\geq p^N.$$
 \end{Proposition3}
 \noindent{\sl Proof} The assertions 1),2),3),4) come from the fact that $\Gamma_{\delta}, \gamma_{-\delta}, U$ are continuous and that these assertions are true for pseudo-polynomials. The assertion 5) is a direct application of Proposition \ref{Proposition2} . $\diamondsuit$\par
 Let's recall the following remarkable result due to W. Sinnott:
 \newtheorem{Proposition4}[Proposition2]{Proposition}
 \begin{Proposition4} \label{Proposition4}
 Let $r_1(T),\cdots ,r_s(T)\in \mathbb F_q(T)\cap \mathbb F_q[[T]].$ Let $c_1,\cdots ,c_s\in \mathbb Z_p\setminus \{ 0\}$ and suppose that:
 $$\sum_{i=1}^s r_i((1+T)^{c_i}-1)=0.$$
 Then:
 $$\forall a\in \mathbb Z_p,\, \sum_{c_i \equiv a\pmod{\mathbb Q^*}}r_i((1+T)^{c_i}-1)\, \in \mathbb F_q.$$
 \end{Proposition4}
 \noindent{\sl Proof} See \cite{SI2}, Proposition 1. $\diamondsuit$\par
 Let's give a first application of this latter result:
 \newtheorem{Proposition5}[Proposition2]{Proposition}
 \begin{Proposition5} \label{Proposition5}
 Let $\delta \in \mathbb Z/(p-1)\mathbb Z$ and let $F(T)\in K(T)\cap \Lambda.$\par
 \noindent 1) If $\delta $is odd or if $\delta =0,$ then:
 $$\mu (\Gamma_{\delta }(F(T)))=\mu (U(F(T))+(-1)^{\delta }U(F((1+T)^{-1}-1))).$$
 2) If $\delta $is even and $\delta \not =0,$ then:
 $$\mu (\Gamma_{\delta }(F(T)))=\mu (U(F(T))+U(F((1+T)^{-1}-1))-2U(F)(0)).$$
 \end{Proposition5}
 \noindent{\sl Proof}  The case $\delta =0$ has already been obtained by Sinnott (\cite{SI1}, Theorem 1). We prove 1), the proof of 2) is quite similar. Now, observe that 1) is a consequence of Proposition \ref{Proposition3} and the following fact:\par
 \noindent Let $F(T)\in K(T)\cap \Lambda ,$ then $\mu(\gamma_{-\delta }(F(T)))=\mu (F(T)+(-1)^{\delta }F((1+T)^{-1}-1)).$\par
 \noindent  Let's prove this fact. Let $r(T)\in \Lambda,$ observe that:
 $$\gamma_{-\delta }(r(T))= (-1)^{\delta} \gamma_{-\delta}(r((1+T)^{-1}-1)).$$
 We can assume that $F(T)+(-1)^{\delta}F((1+T)^{-1}-1)\not =0.$ Write:
 $$F(T)+(-1)^{\delta}F((1+T)^{-1}-1)=\pi^{m}G(T),$$
 where $m\in \mathbb N,$ and  $G(T)\in \Lambda \setminus \pi \Lambda .$ Note that $G(T)\in K(T).$ We must prove that $\gamma_{-\delta}(G(T))\not \equiv 0\pmod{\pi}.$ Suppose that it is not the case, i.e. $\gamma_{-\delta}(G(T)) \equiv 0\pmod{\pi}.$ Then:
 $$G(0)\equiv 0\pmod{\pi}.$$
 Furthermore, by Proposition \ref{Proposition4}, there exists $c\in O_K$ such that:
 $$G(T)+(-1)^{\delta}G((1+T)^{-1}-1)\equiv c\pmod{\pi}.$$
 But, we must have $c\equiv 0\pmod{\pi}.$ Observe that:
 $$G(T)=(-1)^{\delta}G((1+T)^{-1}-1).$$
 Therefore we get $G(T)\equiv 0\pmod{\pi}$ which is a contradiction. $\diamondsuit$\par
  
  \newtheorem{Lemma5}[Proposition2]{Lemma}
  \begin{Lemma5} \label{Lemma5}
  Let $F(T)\in \mathbb F_q(T)\cap \mathbb F_q[[T]].$ Then $F(T)$ is a pseudo-polynomial if and only if there exists some integer $n\geq 0$ such that $(1+T)^n F(T)\in \mathbb F_q[T].$
  \end{Lemma5}
  \noindent{\sl Proof} Assume that $F(T)$ is a pseudo-polynomial. We can suppose that $F(T)\not =0.$ Write:
  $$F(T)=\sum_{i=1}^r c_i (1+T)^{a_i},$$
  where $c_1,\cdots, c_r \in \mathbb F_q^*,$ $a_1,\cdots a_r\in \mathbb Z_p$ and $a_i\not = a_j$ for $i\not =j.$ Since $F(T)\in \mathbb F_q(T)$ there exist $m,n\in \mathbb N\setminus\{0\},$ $m>{\rm Max}\{ v_p(a_i-a_j),\, i\not =j\},$ such that:
  $$(T^{q^n}-T)^{q^{m}}F(T)\in \mathbb F_q[T].$$
  Thus:
  $$\sum_{i=1}^rc_i (1+T)^{a_i+q^{n+m}}-\sum_{i=1}^rc_i (1+T)^{a_i+q^{m}}\, \in \mathbb F_q[T].$$
  Observe that:\par
  \noindent- $\forall i,j \in \{ 1,\cdots ,r\},$ $a_i+q^{n+m}\not = a_j+q^m,$\par
  \noindent - $a_i+q^m=a_j+q^m \Leftrightarrow i=j.$\par
  \noindent Thus, by Lemma \ref{Lemma2}, we get:
  $$\forall i\in \{ 1,\cdots r\}, \, a_i+q^m \in  \mathbb N.$$
  Therefore $(1+T)^{q^m}F(T)\in \mathbb F_q[T].$ The Lemma follows. $\diamondsuit$\par
  Let's give a second application of Proposition \ref{Proposition4}:
  \newtheorem{Proposition6}[Proposition2]{Proposition}
  \begin{Proposition6} \label{Proposition6}
  Let $\delta \in \mathbb Z/(p-1)\mathbb Z$ and let $F(T)\in \mathbb F_q(T)\cap \mathbb F_q[[T]].$ Suppose that there exist an integer $r\in \{ 0,\cdots, (p-3)/2\},$  $c_1,\cdots c_r \in \mathbb Z_p\setminus\{ 0\},$ $G_1(T),\cdots , G_r(T) \in \mathbb F_q(T)\cap \mathbb F_q[[T]]$ and a pseudo-polynomial $R(T)\in \mathbb F_q[[T]]$ such that:
  $$\overline{\gamma_{\delta}}(F(T))=R(T)+\sum_{i=1}^r G_i((1+T)^{c_i}-1).$$
  Then, there exists an integer $n\geq 0$ such that:
  $$(1+T)^n (F(T)+(-1)^{\delta} F((1+T)^{-1}-1))\, \in \mathbb F_q[T].$$
  \end{Proposition6}
  \noindent{\sl Proof} Note that if $\eta,\eta'\in \mu_{p-1}:$ $\eta\equiv \eta'\pmod{\mathbb Q^*}\Leftrightarrow \eta=\eta'\,  {\rm or}\, \eta =-\eta'.$ Since $r<(p-1)/2,$ by Proposition \ref{Proposition4}, there exists $\eta \in \mu_{p-1}$ such that:
  $$\overline{\eta}^{\delta} F((1+T)^{\eta }-1)+\overline{-\eta}^{\delta} F((1+T)^{-\eta }-1)\, {\rm is\, a \, pseudo-polynomial}.$$
  Therefore:
  $$F(T)+(-1)^{\delta} F((1+T)^{-1}-1)\, {\rm is\, a \, pseudo-polynomial}.$$
  It remains to apply Lemma \ref{Lemma5}. $\diamondsuit$\par
 Let $F(T)\in \Lambda.$ We say that $F(T)$ is a pseudo-rational function if $F(T)$ is the quotient of two pseudo-polynomials. For example, $\forall a\in \mathbb Z_p,$ $\forall b\in \mathbb Z_p^*,$ $\frac{(1+T)^a -1}{(1+T)^b-1}$ is a pseudo-rational function. We finish this section by giving a generalization of \cite{SI2}, Theorem1:
 \newtheorem{Theorem1}[Proposition2]{Theorem}
 \begin{Theorem1} \label{Theorem1}
 Let $\delta \in \mathbb Z/(p-1)\mathbb Z$ and let $F(T)\in \mathbb F_q(T)\cap \mathbb F_q[[T]].$ Then $\overline{\Gamma_{\delta}}(F(T))$ is a pseudo-rational function if and only if there exists some integer $n\geq 0$ such that:
 $$(1+T)^n(\overline{U}(F(T))+(-1)^{\delta}\overline{U}(F((1+T)^{-1}-1)))\, \in \mathbb F_q[T].$$
 \end{Theorem1} 
 \noindent{\sl Proof} Assume that $\overline{\Gamma_{\delta}}(F(T))$ is a pseudo-rational function. Then , by 3) of Proposition \ref{Proposition3} and Proposition \ref{Proposition2}, there exist $c_1,\cdots, c_r \in \mathbb F_q^*,$ $a_1,\cdots ,a_r \in \mathbb Z_p,$ $a_i\not = a_j$ for $i\not =j,$ such that:
 $$\overline{\Gamma_{\delta}}\overline{\gamma_{-\delta}}\overline{ U}(\sum_{i=1}^r c_i F((1+T)^{\kappa^{a_i}}-1)))\, {\rm is\, a \, pseudo-polynomial}.$$
 This implies, again by Proposition \ref{Proposition2}, that:
 $$\overline{\gamma_{-\delta}}\overline{ U}(\sum_{i=1}^r c_i F((1+T)^{\kappa^{a_i}}-1)))\, {\rm is\, a \, pseudo-polynomial}.$$
 Set:
 $$G(T)=\overline{U}(F(T))+(-1)^{\delta}\overline{U}(F((1+T)^{-1}-1))\, \in \mathbb F_q(T)\cap \mathbb F_q[[T]].$$
 Now, by Proposition \ref{Proposition4}, there exist $d_1,\cdots ,d_{\ell}\in \mathbb F_q^*,$ $b_1,\cdots b_{\ell}\in \mathbb Z_p,$ $b_i\not =b_j$ for $i\not =j,$ $\eta_1,\cdots , \eta_{\ell}\in \mu_{p-1},$ with $\forall i,j\in \{1,\cdots ,\ell \},$ $\eta_i\kappa^{b_i}\equiv \eta_j\kappa^{b_j}\pmod{\mathbb Q^*},$ and $\eta_i\kappa^{b_i}\not = \eta_j\kappa^{b_j}$ for $i\not = j,$ such that:
 $$\sum_{i=1}^{\ell}d_iG((1+T)^{{\eta_i}\kappa^{b_i}}-1) \, {\rm is\, a \, pseudo-polynomial}.$$
 For $i=1,\cdots , \ell,$ write:
 $$\eta_i\kappa^{b_i}= \eta_1\kappa^{b_1}x_i,$$
 where $x_i\in \mathbb Q^*\cap \mathbb Z_p^*,$ and $x_i\not = x_j$ for $i\not = j.$
 Since $G(T)=(-1)^{\delta} G((1+T)^{-1}-1),$ we can assume that $x_1,\cdots x_{\ell}$ are positives. Now, we get:
 $$\sum_{i=1}^{\ell}d_iG((1+T)^{x_i}-1) \, {\rm is\, a \, pseudo-polynomial}.$$
 Therefore, there exist $N_1,\cdots ,N_{\ell}\in \mathbb N\setminus\{0\},$ $N_i\not =N_j$ for $i\not =j,$ such that:
 $$\sum_{i=1}^{\ell}d_iG((1+T)^{N_i}-1) \, {\rm is\, a \, pseudo-polynomial}.$$
 Now, by Lemma  \ref{Lemma5}, there exists some integer $N\geq 0$ such that:
 $$(1+T)^N(\sum_{i=1}^{\ell}d_iG((1+T)^{N_i}-1)) \, \in \mathbb F_q[T].$$
 But, since $G(T)\in \mathbb F_q(T)\cap \mathbb F_q[[T]],$ $d_1,\cdots ,d_{\ell}\in \mathbb F_q^*,$ $N_1,\cdots N_{\ell }\in \mathbb N\setminus \{ 0\}$ and $N_i\not =N_J$ for $i\not =j,$ this implies that there exist some integer $n\geq 0$ such that $(1+T)^n G(T)\in \mathbb F_q[T].\, \diamondsuit$\par

 %---------------------------------------P-ADIC L-FUNCTIONS--------------------------------------------

 \section{Application to Kubota-Leopoldt $p$-adic L-functions}
 ${}$\par
 Let $\theta$ be a Dirichlet character of the first kind, $\theta \not =1$ and $\theta$ even. We denote by $f(T,\theta)$ the Iwasawa power series attached to the $p$-adic L-function $L_p(s,\theta)$ (see \cite{WAS}, Theorem 7.10). Write:
 $$\theta =\chi \omega^{\delta +1},$$
 where $\chi$ is of conductor $d,$ $d\geq 1$ and $d\not \equiv 0\pmod{p},$ and $\delta \in \mathbb Z/(p-1)\mathbb Z.$ Set $\kappa =1+pd$ and $K=\mathbb Q_p(\chi).$ We set:
 $$F_{\chi}(T)=\frac{\sum_{a=1}^d\chi (a)(1+T)^a}{1-(1+T)^d}.$$
 Let's give the basic properties of $F_{\chi}(T):$
 \newtheorem{Lemma6}{Lemma}[section]
 \begin{Lemma6} \label{Lemma6}
 ${}$\par
 \noindent 1) If $d\geq 2,$ $F_{\chi}(T)\in \Lambda.$\par
 \noindent 2) If $d=1,$ $\forall \alpha \in \mathbb Z/(p-1)\mathbb Z,$ $\alpha \not =1,$ $\gamma_{\alpha}(F_{\chi }(T))\in \Lambda.$\par
 \noindent 3) $U(F_{\chi}(T))=F_{\chi}(T)-\chi(p) F_{\chi}((1+T)^p-1).$\par
 \noindent 4) If $d\geq 2,$ $F_{\chi}((1+T)^{-1}-1)=\varepsilon F_{\chi}(T),$ wher $\varepsilon =1$ if $\chi$ is oddd and $\varepsilon=-1$ if $\chi$ is even.\par
 \noindent 5) If $d=1,$ $F_{\chi}((1+T)^{-1}-1)=-1-F_{\chi}(T).$\par
 \end{Lemma6}
 \noindent {\sl Proof}  1), 4) and 5) are obvious.\par
 \noindent 2)  For $d=1,$ we have:
 $$F_{\chi}(T)=-1+\frac{\sum_{a=0}^{p-1}(1+T)^a}{1-(1+T)^p}.$$
 Set:
 $$G(T)=(1-(1+T)^p)\gamma_{\alpha }(F_{\chi}(T)).$$
 Note that:
 $$\forall \eta \in \mu_{p-1},\, \frac{1-(1+T)^p}{1-(1+T)^{\eta p}}\equiv \eta^{-1} \pmod{\omega_1(T)}.$$
 Therefore:
 $$(p-1) G(T) \equiv \sum_{\eta\in \mu_{p-1}}\eta^{\alpha -1}\sum_{a=0}^{p-1}(1+T)^{\eta a}\pmod{\omega_1(T)}.$$
 Thus:
 $$(p-1) G(T) \equiv \sum_{\eta\in \mu_{p-1}}\eta^{\alpha -1}\sum_{b=0}^{p-1}(1+T)^{b}\pmod{\omega_1(T)}.$$
 Since $\alpha \not =1,$ we get:
 $$G(T)\equiv 0\pmod{\omega_1(T)}.$$
 Therefore $\gamma_{\alpha}(F_{\chi}(T))\in \Lambda.$\par
 \noindent 3) For $d=1,$ we have:
 $$U(F_{\chi}(T))=\frac{\sum_{a=1}^{p-1}(1+T)^a}{1-(1+T)^p}= F_{\chi} (T) -F_{\chi }((1+T)^p-1).$$
 Now, let $d\geq 2.$ Set $q_0=\kappa=1+pd.$ Note that:
 $$F_{\chi}(T)=\frac{\sum_{a=1}^{q_0}\chi(a) (1+T)^a}{1-(1+T)^{q_0}}.$$
 Therefore:
 $$U(F_{\chi}(T))=\frac{\sum_{a=1,\, a\not \equiv 0\pmod{p}}^{q_0}\chi(a) (1+T)^a}{1-(1+T)^{q_0}}.$$
 But:
 $$F_{\chi}(T)-\chi (p)F_{\chi}((1+T)^p-1)= \frac{\sum_{a=1}^{q_0}\chi(a) (1+T)^a}{1-(1+T)^{q_0}}-\chi (p)\frac{\sum_{a=1}^{d}\chi(a) (1+T)^{pa}}{1-(1+T)^{q_0}}.$$
 The Lemma follows easily. $\diamondsuit$\par
 \newtheorem{Lemma7}[Lemma6]{Lemma}
 \begin{Lemma7} \label{Lemma7}
 Assume that $d\geq 2.$ The denominator of $F_{\chi}(T)$ is $\phi_{d}(1+T)$ where $\phi_d(X)$ is the $d$th cyclotomic polynomial and the same is true for $\overline{F_{\chi}(T)}.$\par
 \end{Lemma7}
 \noindent{\sl Proof}  Let $\zeta \in \mu_d.$ If $\zeta$ is not a primite $d$th root of unity, then, by \cite{WAS}, Lemma 4.7, we have:
 $$\sum_{a=1}^d \chi (a) \zeta^a=0.$$
 If $\zeta $ is a primitive $d$th root of unity, then by \cite{WAS}, Lemma 4.8, we have:
 $$\sum_{a=1}^d \chi (a) \zeta ^a \not \equiv 0 \pmod{\widetilde{\pi }},$$
 where $\widetilde {\pi }$ is any prime of $K(\mu_d ).$ $\diamondsuit$\par
 \newtheorem{Lemma8}[Lemma6]{Lemma}
 \begin{Lemma8} \label{Lemma8}
 The derivative of $\gamma_{-\delta}(F_{\chi}(T))$ is not a pseudo-polynomial modulo $\pi.$\par
 \end{Lemma8}
 \noindent{\sl Proof}  We first treat the case $d\geq 2.$ By 3) and 4) of  Lemma \ref{Lemma6}, Lemma \ref{Lemma7} and Proposition \ref{Proposition6}, $\overline{\gamma_{-\delta}}\overline{U} (\overline{F_{\chi}(T)})$  is not a pseudo-polynomial. But observe that $\overline{U}=\overline{D^{p-1}}. $ Thus $\overline{D}\overline{\gamma_{-\delta}} (\overline{F_{\chi}(T)})$ is not a pseudo-polynomial.\par
 \noindent For the case $d=1.$ Set $\widetilde {F_{\chi}(T)}=F_{\chi}(T)-2F_{\chi}((1+T)^2-1)=1-\frac{1}{2+T}.$ Observe that:\par
 \noindent - $ \widetilde {F_{\chi}((1+T)^{-1}-1)}=1-\widetilde {F_{\chi}(T)},$\par
 \noindent -  $U(\widetilde {F_{\chi}(T)})=\widetilde {F_{\chi}(T)}-\widetilde {F_{\chi}((1+T)^p-1)}.$\par
 \noindent Therefore, as in the case $d\geq2,$ $ \overline {\gamma_{-\delta}}\overline{U} (\overline{\widetilde{F_{\chi}(T)}})$  is not a pseudo-polynomial.  Thus $ \overline{\gamma_{-\delta}}\overline{U} (\overline{F_{\chi}(T)})$  is not a pseudo-polynomial.  And one can conclude as in the case $d\geq 2.$ $\diamondsuit$ \par
 \newtheorem{Lemma9}[Lemma6]{Lemma}
 \begin{Lemma9} \label{Lemma9}
 $$\Gamma_{\delta}\gamma_{-\delta}(F_{\chi}(T))=f(\frac{1}{1+T}-1,\theta).$$
 \end{Lemma9}
 \noindent{\sl Proof} We treat the case $d=1,$ the case $d\geq 2$ is quite similar. Set $T=e^Z-1.$ We get:
 $$\gamma_{-\delta}(F_{\chi}(T))=\sum_{n\geq 0,\, n\equiv 1+\delta \pmod{p-1}}\frac{B_n}{n!}Z^{n-1}.$$
 Thus, by \cite{WAS}, Theorem 5.11, we get:
 $$\forall k\in \mathbb N, k\equiv \delta \pmod{p-1},\, D^k\gamma_{-\delta}U(F_{\chi}) (0)=L_p(-k,\theta).$$
 But, by Proposition \ref{Proposition1},we have for $s\in \mathbb Z_p:$
 $$\Gamma_{\delta}\gamma_{-\delta}U(F_{\chi})(\kappa^s-1)={\rm lim}_n D^{k_n(s,\delta)}\gamma_{-\delta}U(F_{\chi})(0)=L_p(-s,\theta)=f(\kappa^{-s}-1,\theta).$$
 The Lemma follows. $\diamondsuit$\par
 We can  now state and prove our main result:
 \newtheorem{Theorem2}[Lemma6]{Theorem}
 \begin{Theorem2} \label{Theorem2}
 ${}$\par
 \noindent 1) $\overline{f(T,\theta )}$ is not a pseudo-rational function.\par
 \noindent 2) $\lambda (f(T,\theta ))< (\frac{p-1}{2}\phi(d))^{\phi (p-1)},$ where $\phi$ is Euler's totient function.\par
 \end{Theorem2}
 \noindent{\sl Proof}${}$\par
 \noindent 1) Assume the contrary, i.e.  $\overline{f(T,\theta )}$ is a pseudo-rational function. Then $\overline{f(\frac{1}{1+T}-1,\theta )}$ is also a pseudo-rational function. Thus $\overline{\Gamma_{\delta}}\overline {\gamma_{-\delta}}\overline{U} (\overline {F_{\chi}(T)})$ is a pseudo-rational function.\par
 \noindent  We first treat the case $d\geq 2.$ By Theorem \ref{Theorem1}, there exists an integer $n\geq 0$ such that $(1+T)^n(\overline{U}(\overline{F_{\chi}(T)})+(-1)^{\delta }\overline{U}(\overline{F_{\chi}((1+T)^{-1}-1)}) \in \mathbb F_q[T].$ This is a contradiction by 3) and 4) of  Lemma \ref{Lemma6} and Lemma \ref{Lemma7}.\par
 \noindent For the case $d=1.$ We work with $\widetilde {F_{\chi}(T)}=F_{\chi}(T)-2F_{\chi}((1+T)^2-1)=1-\frac{1}{2+T}.$ Then , by Proposition \ref{Proposition3}, $\overline{\Gamma_{\delta}}\overline {\gamma_{-\delta}}\overline{U} (\overline {\widetilde{F_{\chi}(T)}})$ is a pseudo-rational function. We get a contradiction as in the case $d\geq 2.$\par
 \noindent 2) Our proof is inspired by a method introduced by S. Rosenberg (\cite{ROS}). We first treat the case $d=1.$ Note that we can assume that $\lambda(f(T,\theta ))\geq 1.$ Now, by Lemma \ref{Lemma8}:
 $$\mu(\gamma_{-\delta}(F_{\chi}(T)))=0.$$
 Futhermore, we have:
 $$\gamma_{-\delta}(F_{\chi})(0)\equiv 0\pmod{\pi }.$$
 Therefore, by 3) of Lemma \ref{Lemma6}, we get:
 $$\lambda(\gamma_{-\delta}U(F_{\chi}(T)))=\lambda (\gamma_{-\delta}(F_{\chi}(T))).$$
 Therefore we have to evaluate $\lambda (\gamma_{-\delta}(F_{\chi}(T))).$ Set $F(T)=\frac{-1}{T}.$ Since $\delta$ is odd, we have:
 $$\gamma_{-\delta}(F_{\chi}(T))=\gamma_{-\delta}(F(T)).$$
 Observe that $F((1+T)^{-1}-1)=1-F(T).$ Let $S\subset \mu_{p-1}$ be  a set of representatives of $\mu_{p-1}/\{ 1,-1\}.$ We have:
 $$(p-1)\gamma_{-\delta}(F(T))=2\sum_{\eta \in S} \eta^{-\delta}F((1+T)^{\eta}-1)\, -\sum_{\eta \in S}\eta^{-\delta}.$$
 Set:
 $$G(T)=(\prod_{\eta\in S}((1+T)^{\eta}-1))\gamma_{-\delta}(F(T)).$$
 Then:\par
 \noindent- $\mu (G(T))=0,$\par
 \noindent- $\lambda (G(T))=\frac{p-1}{2}+\lambda (\gamma_{-\delta }(F(T))).$\par
 \noindent For $S'\subset S,$ write $t(S')=\sum_{x\in S'}x.$ We can write:
 $$G(T)=\sum_{S'\subset S} a_{S'} (1+T)^{t(S')},$$
 where $a_{S'}\in O_K.$ Set:
 $$N={\rm Max}\{v_p(t(S')-t(S'')),\, S',S''\subset S,\, t(S')\not = t(S'')\}.$$
 It is clear that:
 $$p^N<(\frac{p-1}{2})^{\phi (p-1)}.$$
 But, by Lemma \ref{Lemma2}, we have:
 $$\lambda(G(T))<p^{N+1}.$$
 Thus, by Propositon \ref{Proposition3}, we get:
 $$\lambda(f(T,\theta))=\lambda(f(\frac{1}{1+T}-1,\theta))<p^N<(\frac{p-1}{2})^{\phi (p-1)}.$$
 Now, we treat the general case, i.e. $d\geq 2.$ Again we can assume that $\lambda(f(T,\theta ))\geq 1.$ Thus as in the case $d=1,$ we get:
 $$\lambda(\gamma_{-\delta}U(F_{\chi}(T)))=\lambda (\gamma_{-\delta}(F_{\chi}(T))).$$
 Now, by Lemma \ref{Lemma7}, we can write:
 $$F_{\chi}(T)=\frac{\sum_{a=0}^{\phi(d)-1}r_a (1+T)^a}{\phi_d(1+T)},$$
 where $r_a\in O_K$ for $a\in \{ 0,\cdots, \phi(d)-1\}.$
 Let again $S\subset \mu_{p-1}$ be a set of representatives of $\mu_{p-1}/\{ 1,-1\}.$ By Lemma \ref{Lemma6}, we have:
 $$(p-1)\gamma_{-\delta}(F_{\chi}(T))=2\sum_{\eta \in S}\eta^{-\delta}F_{\chi}((1+T)^{\eta }-1).$$
 Set:
 $$G(T)=(\prod_{\eta \in S}\phi_d((1+T)^{\eta})))\gamma_{-\delta}(F_{\chi}(T)).$$
 We have:
 $$G(T)=\sum_{a=0}^{\phi(d)-1}\sum_{\eta \in S}\sum_{S'\subset S\setminus \{ \eta \}}\, \, \sum_{\underline{d}=(d_{\eta'})_{\eta'\in S'},\, d_{\eta'}\in \{ 0,\cdots, \phi (d)\}}b_{S',\underline{d}}(1+T)^{a\eta +\sum_{\eta'\in S'}d_{\eta'}\eta'},$$
 where $b_{S',\underline{d}}\in O_K.$ Note that again $\mu (G(T))=0$ and that $\lambda (G(T))=\lambda (\gamma_{-\delta}(F_{\chi}(T))).$ Now, for $a,b \in \{ 0,\cdots , \phi (d)-1\},$ $\eta_1,\eta_2\in S,$ $S_1\in S\setminus \{ \eta_1\},$ $S_2\in S \setminus\{ \eta_2\},$ set:
 $$V=a\eta_1+\sum_{\eta \in S_1} d_{\eta}\eta \, -b\eta_2-\sum_{\eta \in S_2} d_{\eta }' \eta,$$
 where $\forall \eta \in S_1,$ $d_{\eta}Ê\in \{ 0, \cdots , \phi (d)\},$ and $\forall \eta \in S_2,$ $d_{\eta}'Ê\in \{ 0, \cdots , \phi (d)\}.$\par
 \noindent If $\eta_1=\eta_2$ then we can write:
 $$V=(a-b)\eta_1+\sum_{\eta \in S'} u_{\eta }\eta,$$
 where $\mid u_{\eta}\mid \in \{ 0,\cdots ,\phi(d)\}$ and $\mid S'\mid \leq \frac {p-3}{2}.$\par
 \noindent  If $\eta_1\not =\eta_2,$ we can write:
 $$V=a'\eta_1+b'\eta_2+\sum_{\eta\in S'}u_{\eta} \eta,$$
 where $\mid a'\mid, \mid b'\mid, \mid u_{\eta }\mid \in\{ 0,\cdots ,\phi (d)\},$ and $\mid S'\mid\leq \frac{p-5}{2}.$
 Therefore, if $V\not = 0,$ we get:
 $$p^{v_p(V)}<(\frac{p-1}{2}\phi(d))^{\phi(p-1)}.$$
 Now, we can conclude as in the case $d=1.$ $\diamondsuit$\par
 Let $E$ be a number field and let $E_{\infty}/E$ be the cyclotomic $\mathbb Z_p$-extension of $E.$ For $n\geq 0,$ let $A_n$ be the $p$th Sylow subgroup of the ideal class group of the $n$th layer in $E_{\infty}/E.$Then , by \cite{WAS}, Theorem 13.13, there exist $\mu_p(E)\in \mathbb N, \lambda_p(E)\in \mathbb N$ and $\nu_p(E) \in \mathbb Z,$ such that for all sufficiently large $n:$
 $$\mid A_n\mid =p^{\mu_p(E)p^n+\lambda_p(E)n+\nu_p(E)}.$$
 Recall that it is conjectured that $\mu_p(E)=0$ and if $E$ is an abelian number field it has been proved by B. Ferrero and L. Washington (\cite{FW}).
 \newtheorem{Theorem3}[Lemma6]{Corollary}
 \begin{Theorem3} \label{Theorem3}
 Let $F$ be an abelian number field of conductor $N.$ Write $N=p^m d,$ where $m\in \mathbb N$ and $d\geq 1,$ $d\not \equiv 0\pmod{p}.$ Then:
 $$\lambda_p(F)<2(\frac{p-1}{2}\phi(d))^{\phi(p-1)+1}.$$
 \end{Theorem3}
 \noindent{\sl Proof} Set, for all $n\geq 0,$  $q_n=p^{n+1}d.$ Then $F\subset \mathbb Q(\mu_{q_m}).$ It is not difficult to see that (see the arguments in the proof of Theorem 7.15 in \cite{WAS}) :
 $$\lambda_p(F)\leq \lambda_p(\mathbb Q(\mu_{q_m})).$$
 But, note that:
 $$\lambda_p(\mathbb Q(\mu_{q_m}))=\lambda_p(\mathbb Q(\mu_{q_0})).$$
 Now, by \cite {WAS} Proposition 13.32 and Theorem 7.13:
 $$\lambda_p(\mathbb Q(\mu_{q_0}))\leq 2\sum_{\theta \, {\rm even},\, \theta \not =1,\, f_{\theta}\mid q_0}\lambda (f(T,\theta)).$$
 It remains to apply Theorem \ref{Theorem2}. $\diamondsuit$\par
 Note that the bound of this latter Corollary is certainly far from the truth even in the case $p=3$ (see \cite{KW}).\par
 
 %---------------------------------------REFERENCES------------------------------------------------------

 \end{document}